\documentclass[12pt]{amsart}
\usepackage{ryanmacros}
\input{ryanquivers.sty}
\usepackage[T1]{fontenc}
\usepackage[light,math]{anttor}

\usepackage{tikz}
\usetikzlibrary{matrix,arrows,backgrounds,shapes.misc,shapes.geometric,patterns,calc,positioning}

\usepackage[colorlinks=true, pdfstartview=FitV, linkcolor=blue, citecolor=blue, urlcolor=blue]{hyperref}

\usepackage{fullpage}
\newcommand{\eq}[1]{E_{#1}}
\newcommand{\eg}[1]{E_{\geq#1}}
\newcommand{\mq}[1]{M_{#1}}

\newcommand{\subrep}[3]{\rep(#1 #2 \subset #3)}
\newcommand{\quotrep}[3]{\rep(#1 #2 \onto #3)}
\newcommand{\Gr}[3]{{\rm Gr}_{#2}^{#3}(#1)}
\DeclareMathOperator{\GL}{GL}
\newcommand{\keyw}[1]{\emph{#1}}
\newcommand{\rcol}{\underline{r}}
\newcommand{\dcol}{\underline{d}}
\newcommand{\ranklocus}[3]{RL(#1, #2, #3)}
\newcommand{\maxint}[1]{[#1]_{+}}

\begin{document}
\title{Rank loci in representation spaces of quivers}
\author{Ryan Kinser}
\address{Department of Mathematics, University of Connecticut, Storrs, CT 06269}
\email{ryan.kinser@uconn.edu}
%\thanks{This material was based upon work supported under NSF Grant DMS 0349019}

\begin{abstract}
Rank functors on a quiver $Q$ are certain additive functors from the category of representations of $Q$ to the category of finite-dimensional vector spaces.  Composing with the dimension function on vector spaces gives a rank function on $Q$.  These induce functions on $\rep(Q, \alpha)$, the variety of representations of $Q$ of dimension vector $\alpha$, and thus can be used to define ``rank loci'' in $\rep(Q, \alpha)$ as collections of points satisfying finite lists of linear inequalities of rank functions.
Although quiver rank functions are not generally semicontinuous like the rank of a linear map, we show here that they do have the geometric property that these rank loci are constructible subvarieties.
The same is true for loci defined by rank functions in Schofield's subrepresentation bundles on $\rep(Q, \alpha)$, and in quiver Grassmannians.
\end{abstract}
%\comment{
\maketitle

%%%%%%%%%%%%%%%%%%%%%%%%%%%%%%%%%%%%%%%%%%%%%%%%%%%
%								INTRODUCTION							%
%%%%%%%%%%%%%%%%%%%%%%%%%%%%%%%%%%%%%%%%%%%%%%%%%%%
\section{Introduction}\label{sect:intro}
There is a rich body of work on quiver representations from both algebraic and geometric viewpoints (see articles such as \cite{MR718127,  MR897322, Nakajima:1996ys,Reineke:2008fk}).  The goal of this paper is to establish a geometric property of \keyw{quiver rank functions}, tools which previously have been used to study tensor products and other algebraic aspects of quiver representations.  We work over an arbitrary field $K$ throughout the paper.

A \keyw{quiver} is just another name for a finite directed graph (possibly with loops, parallel edges, etc.) and a \keyw{representation} of a quiver $Q$ is an assignment of a finite-dimensional vector space to each vertex and a linear map to each arrow of $Q$ (Section \ref{sect:background} covers background and establishes notation).
%(between the spaces associated to the head and tail of the arrow).  
Maps between $Q$ and other quivers give rise to associated quiver rank functions on $Q$.
% (details may be found in Section \ref{sect:rank}).  
These are generalizations of the classical rank of a linear map in that they assign to each representation of $Q$ a nonnegative integer which, roughly, measures the dimension of the largest vector space which is ``propagated'' in some way through the representation.
Rank functions are additive with respect to direct sum and certain ones are multiplicative with respect to the pointwise tensor product of representations.
%\todo{depends how i define terminology}
They have been used to study representation rings of quivers; for example, the multiplicative rank functions on a rooted tree quiver are in bijection with a complete set of primitive, orthogonal idempotents in its representation ring \cite{kinserrootedtrees}.

If we consider the space of matrices of a fixed size $m \times n$, allowing the entries to vary in the field $K$, we get an algebraic variety $M_{m,n}$ on which the classical rank function is semicontinuous (with respect to the Zariski topology, which we use throughout).
% on $M_{m,n}$ which gives it a stratification.
In the quiver setting, if we fix a \keyw{dimension vector} for $Q$ (i.e., a non-negative integer for each vertex), we can take matrices of appropriate sizes over each arrow and allow their entries to vary to get every representation of $Q$ with vector spaces of the prescribed dimensions.  This is the \keyw{representation space} of $Q$ of dimension vector $\alpha$, written $\rep(Q,\alpha)$ or $\rep(\alpha)$ (see Section \ref{sect:repspace}).  As an algebraic variety, it is just isomorphic to an affine space, but it carries the action of a base change group whose orbits are in bijection with the isomorphism classes of representations of $Q$ of dimension vector $\alpha$.
Since rank functions for quivers are defined in terms of representation theory (using certain left and right approximation functors), it is not clear that they are geometric in any sense analogous to classical rank.  One can give examples showing that generalized rank functions are not semicontinuous on $\rep(Q,\alpha)$, but in specific cases they can often be described by vanishing and non-vanishing of some collections of matrix minors.  In these examples, $\rank_Q$ will denote the ``global rank function'' of $Q$, which is used to construct other rank functions (Section \ref{sect:rank}).

\begin{example}\label{ex:notsc}
Let $Q$ be the type $A_3$ quiver $\twosubspaceq$.  Then it is straightforward to compute from the definition that
\[
\rank_Q \left(\twosubspacemaps{K^n}{K^m}{K^r}{A}{B} \right) = \dim_K (\im A \cap \im B ) ,
\]
which is not (in general) semi-continuous on representation spaces.  For example, using the dimension vector $(n,m,r) = (1,2,1)$ and representations $X, Y, Z$ given by
\[
(A, B) = \left( \twobyone{1}{0}, \twobyone{0}{1} \right), \qquad \left( \twobyone{1}{0}, \twobyone{1}{0}\right), \qquad \left( \twobyone{1}{0}, \twobyone{0}{0}\right), \qquad \text{respectively},
\]
we find that $\rank_Q(X) = \rank_Q(Z) = 0$ while $\rank_Q(Y) =1$.  But $Z$ is in the orbit closure of $Y$, which in turn is in the orbit closure of $X$, demonstrating that $\rank_Q$ is neither upper- nor lower-semicontinuous on this representation space.
\end{example}

Sometimes there is no simple description of the global rank function in terms of dimensions of a finite number of kernels, images, etc.

\begin{example}\label{eg:doubleloop}
Let $Q$ be the double loop quiver, so a representation is of the form
\[
\begin{tikzpicture}
%\twoloopquivermaps{K^n}{A}{B}
\draw (0,0) node {$K^n$}; 
\draw [->] (0.1,0.2) arc (150:-140:10pt) ;
\draw (1, 0) node {$B$};
\draw [->] (-0.1,0.2) arc (30:300:10pt);
\draw(-1,0) node {$A$};
\end{tikzpicture}
\]
where $A,B$ are $n\times n$ matrices.  Denote by $A_{\neq 0}$ the largest $A$-stable subspace of $K^n$ whose intersection with $\ker A$ is trivial (i.e., the sum of the generalized eigenspaces of $A$ corresponding to nonzero eigenvalues), and denote by $A_0$ the largest subspace of $K^n$ killed by some power of $A$ (the generalized eigenspace of $A$ with eigenvalue 0).  So we have $K^n = A_0 \oplus A_{\neq 0}$, and also $K^n = B_0 \oplus B_{\neq 0}$ similarly.  Then $\rank_Q (V)$ is the dimension of the largest subspace of $A_{\neq 0} \cap B_{\neq 0}$ which is stable under both $A$ and $B$, modulo the smallest subspace of $K^n$ which is stable under $A,B$ and contains $A_0 + B_0$.
\end{example}

Examples like these lead one to guess that quiver rank functions have some geometric behavior at least. Recall that a subset of a variety $X$ is said to be \keyw{constructible} if it can be obtained from a finite number of subsets of $X$, each of which is either open or closed in $X$, via unions and intersections \cite[Ex.~II.3.18]{MR0463157}; a function $f\colon X\to \Z$ is constructible if its image is finite and each subset $\setst{x\in X}{f(x) =n}$ is constructible.  The significance of constructibility is that this property is preserved by images (and inverse images) of regular maps between algebraic varieties. 
Also, for example, the Euler characteristic of a complex algebraic variety is additive with respect to a partition into constructible subvarieties.  That is, if $X$ is a complex algebraic variety and $X=\coprod X_i$ with each $X_i$ constructible in $X$, then $\chi(X) = \sum_i \chi (X_i)$ (where $\chi(Y)$ denotes the topological Euler characteristic of a variety $Y$) \cite[\S 4.5]{MR1234037}.

%\todo{depend on orientation, good for studying geometry of quiver reps}

The main results of this paper are summarized as follows.

\begin{theoremnonum}
Rank functions are constructible on representation spaces of quivers (Theorem \ref{thm:repqa}), Kac's moduli spaces of indecomposables (Corollary \ref{cor:mod}), subrepresentation bundles (Theorem \ref{thm:bundle}), and quiver Grassmannians (Corollary \ref{cor:qgr}).
\end{theoremnonum}

\subsection*{Acknowledgements}  The author is grateful to Arend Bayer and Milena Herring for assistance in proving Lemma \ref{lem:gp}, and Nicolas Poettering for pointing out an error in the original version of Example \ref{eg:doubleloop}.

%%%%%%%%%%%%%%%%%%%%%%%%%%%%%%%%%%%%%%%%%%%%%%%%
%					BACKGROUND									%
%%%%%%%%%%%%%%%%%%%%%%%%%%%%%%%%%%%%%%%%%%%%%%%%
\section{Background}\label{sect:background}
In this section, we establish notation and recall the definitions of quiver rank functions and representation spaces.   Basic algebraic facts about quiver representations used throughout this paper can be found in the book \cite{assemetal}, while the article \cite{MR897322} provides a good introduction to the geometric side.  A representation $\phi$ of a quiver $Q$ consists of a list of vector spaces $(V_x)$ indexed by the vertices of $Q$, and a list of linear maps $(\phi_a)$ indexed by the arrows of $Q$.  The map $\phi_a$ goes from the vector space at the tail of $a$ to the vector space at the head of $a$.  There is an appropriate notion of a morphism between two representations of the same quiver, which gives a category $\repq$ of representations of $Q$.  The reason that we use $\phi$ rather than the more common $V$ to denote a representation is that we will be primarily interested in fixing the spaces $V_x$ while letting the maps $\phi_a$ vary.

\subsection{Quiver rank functions}\label{sect:rank}
We briefly review the construction of the \keyw{global rank function} of a quiver here; more detail and examples can be found in \cite{kinserrank}.  A representation $\phi$ of a quiver $Q$ has a unique largest subrepresentation $\surjrep(\phi)$ in which the map assigned to each arrow is an epimorphism.  Dually, it has a unique largest quotient $\injrep(\phi)$ which has an injective map associated to each arrow.  The image of the composition $\surjrep(\phi) \into \phi \onto \injrep(\phi)$, denoted $\rkf(\phi)$, has an isomorphism over each arrow.  Here and throughout the paper we only work with connected quivers, so that this forces the vector spaces associated to the vertices in $\rkf(\phi)$ to have a common dimension; this nonnegative integer is then defined to be the global rank of $\phi$, written $\rank_Q (\phi)$.  It is fairly easy to verify that $\surjrep$, $\injrep$, and $\rkf$ are functors, and so this number depends only on the isomorphism class of $\phi$ in $\repq$.

To get more invariants of a representation (numbers depending only on the isomorphism class), we employ morphisms between quivers.  These are just maps which send vertices to vertices and arrows to arrows in a manner compatible with the heads and tails of the arrows.
For any morphism of quivers $f \colon Q' \to Q$, there is an associated \keyw{pullback} functor $f^* \colon \rep(Q) \to \rep(Q')$ given on $\psi=(W_x, \psi_a) \in \rep(Q)$ by
\begin{equation}\label{eq:pb}
%(f^*\psi)_x := W_{f(x)} \qquad \text{and} \qquad (f^*\psi)_a := \psi_{f(a)}
f^*\psi :=(W_{f(x)}, \psi_{f(a)})
\end{equation}
for each vertex $x$ and arrow $a$ (see Example \ref{eg:pb} below).  The \keyw{pushforward} $f_* \colon \rep(Q') \to \rep(Q)$ is given on $\phi=(V_x, \phi_a)$ by
\begin{equation}\label{eq:pf}
(%f_* \phi)_x := \bigoplus_{y \in f^{-1}(x)} V_y \qquad \text{and} \qquad (f_* \phi)_a := \sum_{b \in f^{-1}(a)} \phi_b 
f_* \phi := \left(\bigoplus_{y \in f^{-1}(x)} V_y , \sum_{b \in f^{-1}(a)} \phi_b \right)
\end{equation}
(where we consider the maps $\phi_a$ to be defined on the total vector space $\bigoplus_{x} V_x$ by taking $\phi_a (V_x) = 0$ when $x \neq ta$).  It is easy to see that $f^*$ commutes with tensor product while $f_*$ does not in general.

A map $f\colon Q' \to Q$ induces a function $\rank_f$ on $Q$ given by
\begin{equation}\label{eq:pbrank}
\rank_f (\psi) = \rank_{Q'} (f^* \psi) \qquad \text{for }\psi \in \rep(Q)
\end{equation}
and a function $\rank^f$ on $Q'$ via
\begin{equation}\label{eq:pfrank}
\rank^f (\phi) = \rank_Q(f_* \phi) \qquad \text{for }\phi \in \rep(Q') .
\end{equation}
While both $\rank_f$ and $\rank^f$ are additive with respect to direct sum, only $\rank_f$ will be multiplicative with tensor product, in general.  Given a sequence of quivers $Q_1, \dotsc Q_n$ and morphisms of quivers
\begin{equation}\label{eq:quivseq}
Q_1 \xto{f_1} Q_2 \xleftarrow{f_2} \cdots \xto{f_{n-1}} Q_n ,
\end{equation}
we can even chain together pushforwards and pullbacks to get a function
\begin{equation}\label{eq:rankfunctiondef}
{\rank^{f_1}}_{f_2} \cdots^{f_{n-1}} (\phi) = \rank_{Q_n} (f_{n-1*} \cdots f_2^* f_{1*} \phi)
\end{equation}
which will at least be additive.  Note that if we compose two quiver morphisms $Q_1 \xto{f} Q_2 \xto{g} Q_3$, we get $(gf)^* = f^* g^*$ and $(gf)_* = g_* f_*$, so there is no loss of generality in only considering chains (\ref{eq:quivseq}) with alternating directions of morphisms.
\begin{definition}
Any function of the form (\ref{eq:rankfunctiondef}) for some sequence of quiver morphisms (\ref{eq:quivseq})  will be called an \keyw{(additive) rank function} on $Q$.
\end{definition}
\begin{remark}
In the papers \cite{kinserrank,kinserrootedtrees}, the term ``rank function'' is only applied to multiplicative rank functions.  Since the results of this paper are not {\it a priori} related to multiplicativity, we use the term more broadly to avoid introducing new terminology for nonmultiplicative functions and unnecessarily complicating the language throughout.
\end{remark}
%, while those of the form  are \keyw{corank functions}.
%\todo{really need better terminology}

\begin{example}\label{eg:pb}
Let $f \colon Q' \to Q$ be given below where the vertex and arrow labels indicate the map $f$:
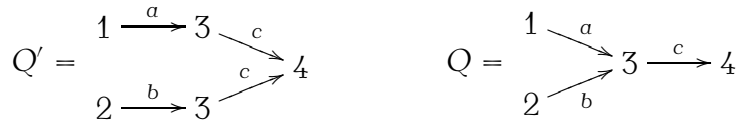
\[
Q'=\QBmaps{1}{2}{3}{3}{4}{a}{b}{c}{c}
\qquad \qquad %\xto{f}\quad
Q=\vcenter{\xymatrix@R=0ex{
{1} \ar[dr]^-{a}  	&	 		& 		\\
			& {3}	\ar[r]^-{c}	& {4}	\\
{2} \ar[ur]_-{b}	&			&		}} 
\]
(e.g., $f$ send both the vertices labeled by 3 in $Q'$ to the one vertex labeled by 3 in $Q$).
Then for $\psi \in \rep(Q)$ we can see the pullback $f^*\psi$ illustrated by
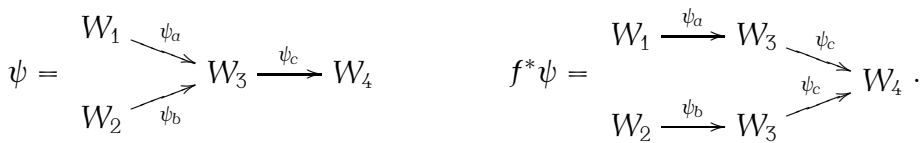
\[
\psi=\vcenter{\xymatrix@R=0ex{
{W_1} \ar[dr]^-{\psi_a} 	&	 		& 		\\
			& {W_3}	\ar[r]^-{\psi_c}	& {W_4}	\\
{W_2} \ar[ur]_-{\psi_b}	&			&		}}
\qquad \qquad
f^*\psi = \QBmaps{W_1}{W_2}{W_3}{W_3}{W_4}{\psi_a}{\psi_b}{\psi_c}{\psi_c} .
\]
In this case, the global rank function of $Q$ can be computed from the definition to be $\rank_Q \psi = \dim \psi_c ( \im \psi_a \cap \im \psi_b) $, while on the other hand $\rank_f \psi = \dim \psi_c \psi_a \cap \psi_c \psi_b$.
\end{example}

\begin{example}
Let $Q'$ be the $n$-subspace quiver and $Q$ of type $A_2$, labeled as
\[
Q'=\nsubspacemaps{1}{2}{n}{0}{a_1}{a_2}{a_n}  \qquad \qquad Q = [n] \xto{a} 0 ,
\]
and $f \colon Q' \to Q$ sending the vertex 0 to 0, and all other vertices to $[n]$.  All the arrows of $Q'$ must collapse to $a$ in $Q$.  The pushforward of $\phi \in \rep(Q')$ can be seen as
\[
\phi=\nsubspacemaps{V_1}{V_2}{V_n}{V_0}{\phi_{a_1}}{\phi_{a_2}}{\phi_{a_n}}  \qquad \qquad f_* \phi = \bigoplus_{i=1}^n V_i \xto{\sum \phi_{a_i}} V_0 ,
\]
and we find that $\rank_{Q'} \phi = \dim \bigcap_i \im \phi_{a_i}$, while $\rank^{f} \phi = \dim \sum_i \im \phi_{a_i}$.
If we first restrict to a subquiver of $Q'$ (pullback along the inclusion), then pushforward along (the restriction of) $f$, we get the functions $\dim \sum_{j \in J} \im \phi_{a_j}$ for any subset $J \subseteq \{1, \dotsc, n\}$.
\end{example}

\subsection{Representation spaces}\label{sect:repspace}
We start by recalling the definitions of $\rep(Q, \alpha)$, the associated base change group, and quiver Grassmannians.
Fix an arbitrary quiver $Q$ and a dimension vector $\alpha$ for $Q$.  Since we will only be interested in a fixed quiver, we often omit $Q$ from the notation.  For an arrow $a$, we let $ta$ and $ha$ be the tail and head of $a$, respectively; for a vertex $x$, denote by $\alpha(x) \in \Z_{\geq 0}$ the component of $\alpha$ at the vertex $x$.  The \keyw{representation space} of $Q$ of dimension vector $\alpha$, written $\rep(Q,\alpha)$ or simply $\rep(\alpha)$, can be defined as
\[
\rep(\alpha) = \bigoplus_{\text{arrows }a} \Hom_K (K^{\alpha(ta)}, K^{\alpha(ha)}) ,
\]
which carries an induced action of the base change group 
\[
\GL(\alpha)= \prod_{\text{vertices }x} \GL_{\alpha(x)}(K) .
\]
A point $\phi \in \rep(\alpha)$ is given by a collection of maps $(\phi_a \colon K^{\alpha(ta)} \to K^{\alpha(ha)})$ indexed by the arrows of $Q$, and two points correspond to isomorphic objects in $\repq$ if and only if they lie in the same orbit of $\GL(\alpha)$.

If $\beta$ is another dimension vector for $Q$, with $\beta(x) \leq \alpha(x)$ for each vertex $x$ (written  $\beta \leq \alpha$), we let
\[
\Gr{\alpha}{\beta}{} = \prod_{\text{vertices }x} \Gr{\alpha(x)}{\beta(x)}{},
\]
where $\Gr{n}{r}{}$ is the classical Grassmannian of $r$-dimensional subspaces of $K^n$.  Thus a point $W \in \Gr{\alpha}{\beta}{}$ is given by a collection of subspaces $(W_x \subseteq K^{\alpha(x)})$.
%(The somewhat untraditional use of $\phi$ to denote points of $\rep(\alpha)$ is meant to make the collections of linear maps more readily distinguishable from collections of subspaces $W$.)
 Then the \keyw{bundle of $\beta$-dimensional subrepresentations} on $\rep(\alpha)$ is the incidence locus
\[
\subrep{}{\beta}{\alpha} = \setst{ (W, \phi) \in \Gr{\alpha}{\beta}{} \times \rep(\alpha)}{\phi_a(W_{ta}) \subseteq W_{ha} \text{ for all arrows }a} .
\]
This construction was introduced by Schofield in \cite[\S3]{schofieldgeneralreps}, where he notes that $\subrep{}{\beta}{\alpha}$ has a projective morphism
\[
p \colon \subrep{}{\beta}{\alpha} \to \rep(\alpha) ,
\]
and is a vector bundle over the homogeneous $\GL(\alpha)$-space $\Gr{\alpha}{\beta}{}$,
\[
q \colon \subrep{}{\beta}{\alpha} \to \Gr{\alpha}{\beta}{} .
\]
The fiber over a representation $\phi \in \rep(\alpha)$ is a projective variety that parametrizes the $\beta$-dimensional subrepresentations of $\phi$, and the fiber over a collection of subspaces $W \in \Gr{\alpha}{\beta}{}$ parametrizes the $\alpha$-dimensional representations which stabilize $W$.

Dually, one can define $\Gr{\alpha}{}{\beta}$ using Grassmannians of quotient spaces $\Gr{n}{}{r}$, and  construct the bundle of $\beta$-dimensional of quotient representations
\[
%\draw (0, 1) node {$\quotrep{}{\alpha}{\beta}$};
%\draw (-2, 0) node {$\Gr{\alpha}{}{\beta}$};
%\draw (2,0) node {$\rep(\alpha)$};
\begin{tikzpicture}[description/.style={fill=white,inner sep=2pt}] 
\matrix (m) [matrix of math nodes, row sep=3em, 
column sep=2.5em, text height=1.5ex, text depth=0.25ex] 
{  & \quotrep{}{\alpha}{\beta} &  \\ 
\Gr{\alpha}{}{\beta} &  & \rep(\alpha) \\ }; 
\path[->,font=\scriptsize] 
(m-1-2) edge node[auto,swap] {$ q' $} (m-2-1) 
%edge node[description] {$ \Psi $} (m-2-3) 
(m-1-2) edge node[auto] {$ p' $} (m-2-3); 
\end{tikzpicture} 
%\xymatrix{ & \quotrep{}{\alpha}{\beta} \ar[dl]^{q'} \ar[dr]_{p'} \\
%\Gr{\alpha}{}{\beta} & & \rep(\alpha) }
\]
with $p'$ projective and $q'$ a vector bundle.

%\subsection{}
For a map between quivers $f \colon Q' \to Q$, the pullback and pushforward functors induce maps between representation spaces of the appropriate dimensions.  From the definitions (\ref{eq:pb}) and (\ref{eq:pf}) we see that these are regular maps of algebraic varieties, so the images of constructible sets under these maps are constructible.
\comment{
Now let $\rcol = (r_i)_{i=1}^n$ be a finite collection of rank and corank functions on $Q$; that is, for each $i$ there exists a quiver $Q_i$ and either a map $f_i \colon Q_i \to Q$ satisfying 
\[
r_i = \rank_{f_i} \colon \repq \to \Z
\]
or a map $f_i \colon Q \to Q_i$ such that 
\[
r_i = \rank^{f_i} \colon \repq \to \Z .
\]
Then any list of nonnegative integers $\dcol = (d_i)_{i=1}^n$ gives a locus in $\rep(\alpha)$ where the specified rank functions $\rcol$ take the values $\dcol$.  Hence, we define the \keyw{rank locus} in $\rep(\alpha)$ with respect to $\rcol$ and $\dcol$ to be
\begin{equation}\label{eq:ranklocus}
\ranklocus{\alpha}{\rcol}{\dcol} = \setst{\phi \in \rep(\alpha)}{r_i (\phi) = d_i \ \text{for each }i}
\end{equation}
\todo{define more generally with less notation}
}
In general, simply looking at all representations where a rank function takes some fixed value may not be very interesting.  So we consider more general loci described by rank functions.
\begin{definition}
A \keyw{rank locus} in $\rep(\alpha)$ is a collection of points satisfying some finite list of linear inequalities in the values of rank functions.
\end{definition}
For a fixed $\alpha$, any rank function on $\rep(\alpha)$ is bounded above by a constant depending on $\alpha$ and the sequence of maps (\ref{eq:quivseq}) used to construct the rank function.
%$\max \{\alpha(x)\}$.  
So from the remarks in the preceding paragraph, we see that rank loci are constructible in general if and only if the global rank function of any quiver is constructible in general.

\begin{example}\label{eg:typea}
When $Q$ is of type $A$, that is, the underlying graph is of the form
\[
\begin{tikzpicture}[point/.style={shape=circle,fill=black,scale=.5pt,outer sep=3pt},>=latex]
   \node[point,label={below:$1$}] (1) at (0,0) {};
   \node[draw, color=white,scale=.6pt,outer sep=3pt] (2) at (2,0) {};
  \node[point,label={below:$2$}] at (2,0) {};
  \node[draw, color=white,scale=.6pt,outer sep=3pt] (3) at (4,0) {};
  \node[point,label={below:$3$}] at (4,0) {};
  \node[point,label={below:$n-1$}] (n-1) at (6,0) {};
%  \node[draw, color=white,scale=.6pt,outer sep=3pt] (3) at (4,0) {};
  \node[point,label={below:$n$}] (n) at (8,0) {};
  
  \path[loosely dotted, line width=1pt]
  	(n-1) edge (3);
  \path[-]
  	(3) edge (2)
  	(2) edge (1);
\path[-]
	(n-1) edge (n);
   \end{tikzpicture}
\]
with any orientation of the arrows, we have that multiplicative rank functions are in bijection with the isomorphism classes of indecomposables.  An isomorphism class in $\repq$ (equivalently, a $\GL(\alpha)$ orbit in $\rep(\alpha)$) is completely determined by the values of these rank functions, so any $\GL(\alpha)$-stable subvariety of $\rep(\alpha)$ can be described as a rank locus.

More specifically, a connected subquiver of $Q$ can be specified by giving its extremal vertices $i$ and $j$, with $1 \leq i \leq j \leq n$. We get a rank function $r_{i,j}$ on $Q$ by restriction to this subquiver (a special case of pullback) then applying the global rank function of the subquiver.  The indecomposable representations $V_{kl}$ of $Q$ are also in bijection with pairs $1 \leq k \leq l \leq n$, and we have that
\[
r_{i,j}(V_{kl}) = 
\begin{cases}
1 & \text{when } k \leq i \leq j \leq l \\
0 & \text{otherwise.}\\
\end{cases}
\]
(This follows from \cite[Theorem~30]{kinserrank}, for example.)  By inclusion-exclusion, we find that the multiplicity of $V_{kl}$ in a representation $V$ is then
\[
r_{k, l}(V)+ r_{k-1,l+1}(V) - r_{k-1, l}(V) - r_{k,l+1}(V)
\]
(where we take $r_{i,j} = 0$ if $i$ or $j$ lie outside $\{1,\dotsc,n\}$), which allows any orbit to be described by rank functions.
\end{example}

\begin{remark}
Example \ref{eg:typea} generalizes to other Dynkin quivers with a ``rooted'' orientation at a minuscule node (see the end of Section 3.4 of \cite{Kinser:2009zr} for a more detailed account).  By also utilizing nonmultiplicative rank functions, the author expects this to work for any Dynkin quiver.  But how to explicitly describe multiplicities of indecomposables with rank functions for a general Dynkin quiver remains an open question.  We also note that while multiplicative rank functions provide information about the tensor product of representations, they might not be the best way of describing rank loci in representation spaces.  For example, Abeasis and Del Fra used certain nonmultiplicative ``rank parameters'' 
(which can be described in terms of our rank functions)
to parametrize the orbits for equioriented $D$ type quivers
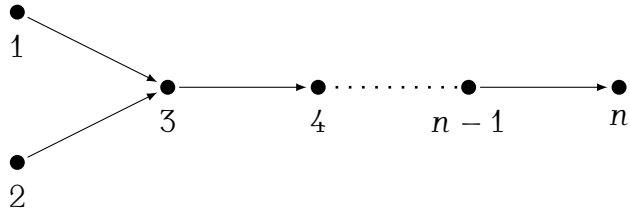
\[
\begin{tikzpicture}[point/.style={shape=circle,fill=black,scale=.5pt,outer sep=3pt},>=latex]
   \node[point,label={below:$1$}] (1) at (0,1) {};
   \node[point,label={below:$2$}] (2) at (0,-1) {};
%      \node[draw, color=white,scale=.6pt,outer sep=3pt] (3) at (2,0) {};
  \node[point,label={below:$3$}] (3) at (2,0) {};
%  \node[draw, color=white,scale=.6pt,outer sep=3pt] (4) at (4,0) {};
  \node[point,label={below:$4$}] (4) at (4,0) {};
  \node[point,label={below:$n-1$}] (n-1) at (6,0) {};
%  \node[draw, color=white,scale=.6pt,outer sep=3pt] (3) at (4,0) {};
  \node[point,label={below:$n$}] (n) at (8,0) {};
  
  \path[loosely dotted, line width=1pt]
  	(n-1) edge (4);
  \path[->]
  	(1) edge (3)
  	(2) edge (3)
	(3) edge (4);
\path[->]
	(n-1) edge (n);
   \end{tikzpicture}
\]
(notice that it is rooted at a minuscule node).  %the orientation is rooted at the minuscule node farthest from the branch point), 
Their functions have the advantage of allowing one to describe degenerations (containment of orbit closures) very easily.

\end{remark}

%%%%%%%%%%%%%%%%%%%%%%%%%%%%%%%%%%%%%%%%%%%%%%%%
%					CONSTRUCT ON REP SPACES							%
%%%%%%%%%%%%%%%%%%%%%%%%%%%%%%%%%%%%%%%%%%%%%%%%
\section{Rank functions on representation spaces}\label{sect:rankrepspace}
Our first goal is to show that the global rank function
\[
\rank_Q \colon \rep(\alpha) \to \N
\]
is constructible.  We will need an intuitive lemma at several points, which we dispose of here.
\begin{lemma}\label{lem:gp}
Let $G$ be an algebraic group, $H$ a closed subgroup, and $p\colon E \to G/H$ a $G$-equivariant vector bundle.  Then for any $G$-equivariant subset $C \subseteq E$, we have that $C$ is closed in $E$ if and only if its intersection with each fiber $F = p^{-1}(gH)$ is closed in $F$.
\end{lemma}
\begin{proof}
Consider the map
\[
\varphi \colon G \times F \to E \qquad (g, f) \mapsto gf .
\]
First, we claim that $C$ is closed in $E$ if and only if $\varphi^{-1}(C)$ is closed in $G \times F$.  
%The ``only if''  implication follows from the continuity of the action, and so for the converse we assume that $\phi^{-1}(C)$ is closed.  
Over any open set $U \subseteq G/H$ which locally trivializes $E$,
%, that is, such that for each $U$ in this cover we have $p^{-1}(U) \simeq U \times F$ with 
we get a diagram
\[
\xymatrix{
\pi^{-1}(U)\times F \ar[r]^{\pi \times id} \ar[d]	& U \times F \ar[d]^{p} \\
\pi^{-1}(U) \ar[r]^{\pi}			& U }
\]
in which we write $\pi\colon G \to G/H$ for the quotient map.
To prove our claim, it is enough to show that $C \cap (U \times F)$ is closed in $U \times F$ if and only if $\varphi^{-1}(C) \cap (\pi^{-1}(U) \times F)$ is closed in $\pi^{-1}(U) \times F$ for any such $U$.  But the $G$-equivariance of $C$ (and the fact that $G$ acts transitively on $G/H$) gives that $C \cap (U \times F) = U \times (C \cap F)$ and $\varphi^{-1}(C) \cap (\pi^{-1}(U) \times F) = \pi^{-1}(U) \times (C \cap F)$, so the claim is verified.

The intersection of a closed subset of $E$ with $F$ is of course closed in $F$.  Now we consider the other projection $\psi \colon G \times F \to F$.  If $C\cap F$ is closed in $F$, then $\psi^{-1}(C \cap F)= G \times (C \cap F)$ is closed in $G \times F$, so $\varphi(G \times (C \cap F)) = C$ is closed in $E$.
\end{proof}

We will be interested in the case where $G= \GL(\alpha)$, $E=\rep(\alpha)$, and $G/H = \Gr{\alpha}{\beta}{}$.  Denote by $\dimv \phi$ the dimension vector of a representation $\phi$ of $Q$, and recall the functors $\surjrep, \injrep$ from Section \ref{sect:intro}.

\begin{definition}
%Let $Q$ be a quiver and fix a dimension vector $\alpha$ for $Q$.  
For each dimension vector $\beta$, we define subsets of $\rep(\alpha)$:
\[
\eq{\beta} =\setst{\phi}{\dimv \surjrep(\phi) = \beta} \qquad \mq{\beta} = \setst{\phi}{\dimv \injrep(\phi) = \beta} .
\]
\end{definition}
These are empty unless $\beta \leq \alpha$.  Say that a representation is \keyw{epimorphic} if each map in it is an epimorphism, so $\surjrep(\phi)$ is the unique maximal epimorphic subrepresentation of $\phi$.

\begin{prop}\label{prop:eqconstructible}
%Let $Q$ be an arbitrary quiver and fix a dimension vector $\alpha$.  Then 
The sets $\eq{\beta}$ and $\mq{\beta}$ are constructible in $\rep(\alpha)$, for any dimension vector $\beta$.
\end{prop}
\begin{proof}
First we will see that the set
\begin{equation}\label{def:X}
X := \setst{(W, \phi)}{\text{the restriction of $\phi$ to $W$ is an epimorphic}}
\end{equation}
is open in $\subrep{}{\beta}{\alpha}$.  Let $F$ be a fiber $q^{-1}(\widetilde{W})$ for some $\widetilde{W} \in \Gr{\alpha}{\beta}{}$.  The intersection $U := X \cap F$ is open in $F$ since it is the locus where the maps given by $\phi \in \rep(\alpha)$ have full rank when restricted to $\widetilde{W}$.
By applying Lemma~\ref{lem:gp} to the complement of $X$, we see that it is open.

The projection $p(X) \subseteq \rep(\alpha)$ is then the set of representations which have some epimorphic subrepresentation of dimension vector $\beta$, and so such a representation has maximal epimorphic subrepresentation of dimension vector at least $\beta$.  So we define
\[
p(X) = \setst{\phi}{\dimv \surjrep(\phi) \geq \beta} =: \eg{\beta}.
\]
By Chevalley's theorem that the images of regular maps of varieties are constructible \cite[Ex.~II.3.19]{MR0463157}, we get that each $\eg{\beta}$ is constructible.   Then it follows that
\[
\eq{\beta} = \eg{\beta} \setminus \bigcup_{\alpha \geq \gamma \gneq \beta} \eg{\gamma}
\]
is constructible also, since the union on the right hand side is finite.

A similar argument utilizing $\quotrep{}{\alpha}{\beta}$
shows that $\mq{\beta}$ is also constructible.
\end{proof}

%\todo{does the next lemma just boil down to some simple statement about sections of projective morphisms?}

\begin{lemma}\label{lem:continsect}
The map 
\begin{align*}
s \colon \eq{\beta} &\to \subrep{}{\beta}{\alpha}  \\
\phi &\mapsto (\surjrep(\phi), \phi)
\end{align*}
is a continuous section of $p$ over $\eq{\beta}$.

Similarly, we have that $s' \colon \mq{\beta} \to \quotrep{}{\alpha}{\beta}$ given by $s'(\phi) = (\injrep(\phi), \phi)$ is a continuous section of $p'$.
\end{lemma}
\begin{proof}
It is clear that $p \circ s$ is the identity on $\eq{\beta}$, so we just need to show that $s$ is continuous.  Retaining the definition of $X$ from (\ref{def:X}) in the proof of Proposition \ref{prop:eqconstructible}, we set $E:= p^{-1}(\eq{\beta})$ for brevity and let $Z := X \cap E$ in $\subrep{}{\beta}{\alpha}$.  That $\im s$ is contained in $Z$ is immediate from the definitions,
and we claim that $Z = \im s$.  A point of $Z$ just a pair $(W, \phi)$ with $W$ a $\beta$-dimensional epimorphic subrepresentation of $\phi$, but such that the unique \emph{maximal} epimorphic subrepresentation $\surjrep(\phi)$ of $\phi$ has dimension $\beta$.  So $W = \surjrep(\phi)$ for such a point, showing that $(W, \phi) \in \im s$.
Thus, $p$ and $s$ give inverse bijections 
\[
\xymatrix{Z \ar@<.5ex>[r]^-{p} & {\eq{\beta}} \ar@<0.5ex>[l]^-{s}} .
\]

The locus $Z$ is open in $E$, since $X$ is open, but we will see that $Z$ is also closed in $E$.  Fixing a collection of subspaces $\widetilde{W} \in \Gr{\alpha}{\beta}{}$, let $F=q^{-1}(\widetilde{W})$ be the fiber over $\widetilde{W}$,
so by Lemma \ref{lem:gp} it is enough to show that $Z_{\widetilde{W}} := Z \cap F$ is closed in $E_{\widetilde{W}} := E \cap F$.  We will do this by constructing it from an intersection of finite unions of closed sets.
Fixing some other $W \neq \widetilde{W} \in \Gr{\alpha}{\beta}{}$, and an arrow $a \in Q\arrows$, we wish to consider the locus in $F$ consisting of pairs $(\widetilde{W}, \phi)$ such that $W$ is a subrepresentation of $\phi$, but $\phi_a$ is not surjective when restricted to $W$.  This is the set
\[
Y(W,a) := \setst{(\widetilde{W}, \phi) \in q^{-1}(\widetilde{W})}{ W \in q(p^{-1}(\phi)) \text{ and }\rank_Q (\phi_a|_{W_{ta}}) < \dim W_{ha} = \beta_{ha}},
\]
which is closed in the vector space $q^{-1}(\widetilde{W})$ because it is given by the vanishing of minors of $\phi_a | _{W_{ta}}$.  Then also the finite union 
\[
Y(W) := \bigcup_{a \in Q\arrows} Y(W, a)
\]
is closed in $q^{-1}(\widetilde{W})$, which can be described as the locus of representations $\phi$ in the fiber over $\widetilde{W}$ which have $W$ as a non-epimorphic subrepresentation.

Now we claim that
\[
Z_{\widetilde{W}} = \bigcap_{W \in q(p^{-1}(\phi)) \setminus \{\widetilde{W}\}} Y(W) \cap E_{\widetilde{W}} ,
\]
which will demonstrate that $Z_{\widetilde{W}}$ is closed in $E_{\widetilde{W}}$.\\
$\subseteq:$ If $(\widetilde{W}, \phi) \in Z_{\widetilde{W}}$, then $\surjrep(\phi) = \widetilde{W}$, so certainly $(\widetilde{W}, \phi) \in E_{\widetilde{W}}$.  For each $W\in q(p^{-1}(\phi)) \setminus \{\widetilde{W}\}$, it is not possible for $W$ to be an epimorphic subrepresentation of $\phi$ because then $W + \widetilde{W} \supsetneq \widetilde{W}$ would be a larger epimorphic subrepresentation, contradicting $\surjrep(\phi) = \widetilde{W}$.  So $\phi_a|_{W_{ta}}$ is not surjective for some arrow $a$, and thus $(\widetilde{W}, \phi) \in Y(W, a) \subseteq Y(W)$.\\
$\supseteq:$  If $(\widetilde{W}, \phi)$ is an element of the right hand side, then in particular it is in $E$ so $\dimv\surjrep(\phi) = \beta$.  But being an element of this intersection says exactly that no other $\beta$-dimensional subrepresentation $W$ is epimorphic, which forces $\surjrep(\phi) = \widetilde{W}$, and so $(\widetilde{W}, \phi) \in Z_{\widetilde{W}}$.

Now we know that $Z$ is closed in $E$.  Since $p$ is a projective morphism, the map $p|_E \colon E \to \eq{\beta}$ obtained by base change is a closed map.  Then it restricts to a closed map on the closed subset $Z$, where it is bijective from above, and thus its inverse $s$ is continuous.
\end{proof}

\begin{example}
Let $Q$ be the loop quiver and consider the dimension vectors $\alpha =2$, $\beta =1$.  Then a point of $\rep(\alpha)$ is given by a $2 \times 2$ matrix, and $\eq{\beta}$ is the locus of matrices which are conjugate to
\[
\begin{pmatrix}
\lambda & 0 \\
0 & 0
\end{pmatrix} , \qquad \lambda \neq 0 .
\]
The fiber of $\subrep{}{\beta}{\alpha} \xto{p} \rep(\alpha)$ over a matrix $M \in \eq{\beta}$ is two points, corresponding to the eigenspaces of $M$, and the bundle $\subrep{}{\beta}{\alpha}$ restricted to $\eq{\beta}$ is isomorphic to two disjoint copies of $\eq{\beta}$.  The section $s$ associates to a matrix the eigenspace with eigenvalue $\lambda$.
\end{example}

Now we are ready to prove the main result.

\begin{theorem}\label{thm:repqa}
For any nonnegative integer $n$, the rank locus
\[
R_n := \setst{\phi \in \rep(\alpha) }{\rank_Q (\phi) = n}
\]
is constructible in $\rep(\alpha)$. Thus, any rank locus
%, as defined in (\ref{eq:ranklocus}), 
is constructible.
\end{theorem}
\begin{proof}
Using Proposition \ref{prop:eqconstructible}, we have a finite partition
\[
\rep(\alpha) = \coprod_{\beta, \gamma \leq \alpha} \eq{\beta} \cap \mq{\gamma}
\]
into constructible sets, so it is enough to show that the intersection of $R_n$ with each set on the right hand side is constructible. From the construction of $\rank_Q$ in Section \ref{sect:rank}, we see that for an arbitrary vertex $x$ the value $\rank_Q (\phi) = \dim_K \rkf (\phi)_x$ is equal to the rank of the linear map 
\[
\surjrep(\phi)_x \into K^{\alpha(x)} \onto \injrep(\phi)_x .
\]
Using Lemma \ref{lem:continsect} and its dual, the composition
\[
\eq{\beta} \cap \mq{\gamma} \xto{s \times s'} \subrep{}{\beta}{\alpha}  \times \quotrep{}{\alpha}{\gamma} 
\xto{q \times q'} \Gr{\alpha}{\beta}{} \times \Gr{\alpha}{}{\gamma}
\]
is continuous, sending $\phi$ to $(\surjrep(\phi), \injrep(\phi))$. Then projecting to the spaces associated to a particular vertex $x$, we get a continuous map 
\begin{align*}
\psi \colon \eq{\beta} \cap \mq{\gamma} &\to \Gr{K^m}{i}{} \times \Gr{K^m}{}{j} \\
\phi &\mapsto (\surjrep(\phi)_x, \injrep(\phi)_x)
\end{align*}
where $m=\alpha(x),\  i=\beta(x)$ and $j = \gamma(x)$.

Now consider the subset 
\[
T_n := 
%\setst{(A,B) \in \Gr{K^m}{i}{} \times \Gr{K^m}{}{j}}{\dim \frac{A+B}{B} = n} = 
\setst{(A, B) \in \Gr{K^m}{i}{} \times \Gr{K^m}{}{j}}{\rank(A \into K^m \onto B) = n}
\]
which is constructible.
Then $R_n \cap E_\beta \cap M_\gamma = \psi^{-1} (T_n)$, and is thus constructible.
\end{proof}

In the representation space $\rep(\alpha)$, many points correspond to isomorphic representations, and in fact the isomorphism classes of representations of $Q$ of dimension vector $\alpha$ are naturally in bijection with the $\GL(\alpha)$ orbits on $\rep(\alpha)$.  So if one wishes to construct a geometric space in which points parametrize some subset of the isomorphism classes of representations of a fixed dimension (i.e., a \keyw{moduli space} of representations), this amounts to putting a geometric structure on some set of orbits in $\rep(\alpha)$.
%Well-known examples of such moduli spaces are those of (semi-)stable representations introduced by King \cite{Kmodulireps} and studied extensively over the past fifteen years \cite{Reineke:2008fk}, and 
One example is the moduli space of indecomposables of a fixed dimension introduced by Kac \cite{MR718127} (via repeated application of Rosenlicht's theorem) and studied in \cite{LeBruyn1986d}.
Since rank functions are constant on orbits, they give well-defined functions on these moduli spaces.
%\todo{false, rank functions to be constant on fibers of quotient map.  Maybe not true for semi-stable case, but probably so for stable case since fibers are closed orbits.  nakajima quiver varieties?!?!?!}

\begin{corollary}\label{cor:mod}
Generalized rank functions are constructible on moduli spaces of indecomposable representations.
\end{corollary}
\begin{proof}
If $M$ is a connected component of a moduli space of indecomposable representations of dimension $\alpha$, then there is a constructible set $U \subset \rep(\alpha)$ and a surjective regular map $U \to M$.  A rank locus in $M$ is the image under this map of the intersection of $U$ and a rank locus in $\rep(\alpha)$.  By Chevalley's theorem, it is then constructible.
\end{proof}

\section{Rank functions on quiver Grassmannians}
A point $(W, \phi) \in \subrep{}{\alpha}{\beta}$ gives a representation $\phi |_W = \left( W_x, \phi_a |_{W_{ta}}\right)$ of $Q$ by taking the vector spaces to be the spaces $(W_x)$, and the maps to be the restrictions of $(\phi_a)$ to these spaces.  Any rank function on $Q$ may be applied to such a point, so we get induced rank functions on $\subrep{}{\alpha}{\beta}$. Each point also gives a quotient representation of $\phi$ simply by modding out the subrepresentation we just considered, and we can apply rank functions to this quotient.  Now by considering inequalities among rank functions, we get rank loci in $\subrep{}{\alpha}{\beta}$ just like we did for $\rep(\alpha)$.
We show that these loci are constructible also.

\begin{theorem}\label{thm:bundle}
Rank loci are constructible subvarieties of the bundle $\subrep{}{\alpha}{\beta}$.
\end{theorem}
\begin{proof}
Let $\pi \colon \subrep{}{\beta}{\alpha} \to \Gr{\alpha}{\beta}{}$ be the projection and consider a fiber $F = \pi^{-1} (\widetilde{W})$. 
% = \setst{(\widetilde{W}, \phi)}{\phi \text{ stabilizes } \widetilde{W}}.
An element $(\phi, \widetilde{W})$ of this fiber gives a representation $\phi |_{\widetilde{W}}$, so we get a regular map
 \[
 b \colon F \to \rep(\beta) %\qquad (\widetilde{W}, V) \mapsto V|_{\widetilde{W}} .
 \] 
 in this way.
%This map is surjective since any representation of $\beta$ can be obtained via extension by 0 to a complement of $\widetilde{W}$. %don't need this i dont' think....
The intersection of a rank locus in $\subrep{}{\beta}{\alpha}$ with $F$ is by definition precisely the preimage of a rank locus of $\rep(\beta)$ under $b$, and thus constructible in $F$.  Furthermore, since $\GL(\alpha)$ acts transitively on the base space $\Gr{\beta}{\alpha}{}$, we find that the rank loci in $\subrep{}{\beta}{\alpha}$ are the $\GL(\alpha)$-orbits of their intersections with $F$.  By writing a rank locus in $F$ using unions and intersections of some open subsets and some closed subsets of $F$, we may apply Lemma \ref{lem:gp} to smear these around and see that the rank locus in $\subrep{}{\beta}{\alpha}$ can be written in the same way.
\end{proof}

A fiber $p^{-1}(\phi)$ is known as a \keyw{quiver Grassmannian}, written $\Gr{\phi}{\beta}{}$.  The isomorphism class of this variety only depends on the isomorphism class of $\phi$, so sometimes we write $\Gr{V}{}{\beta}$ for $V \in \repq$.
These are important in the study of cluster algebras  \cite{Derksen:2009kl,Kclusterquivertriang}.
%\todo{what else? refs?}
Since a rank locus in a quiver Grassmannians is just the intersection of a rank locus in $\subrep{}{\beta}{\alpha}$ with a (closed) fiber of $p$, we get the following corollary.

\begin{corollary}\label{cor:qgr}
Rank loci in quiver Grassmannians are constructible.
\end{corollary}

\begin{example}\label{eg:strat}
Let $Q$ be the Kronecker quiver 
\[
Q = \krontwomaps{1}{2}{a}{b}
\]
and denote by $P_n$ the indecomposable preprojective representation of dimension $(n, n+1)$.  Similarly we write $I_n$ for the indecomposable preinjective of dimension $(n+1, n)$, and $R_n$ for the indecomposable representation of dimension $(n,n)$ in which the map over the bottom arrow is not an isomorphism (it is given by a single Jordan block of eigenvalue 0).
Any submodule of $U \subseteq R_n$ is isomorphic to one of the form $P \oplus R_{k(U)}$ with $P$ a direct sum of preprojective indecomposables and $0 \leq k(U) \leq n$ (of course, $P$ also depends on $U$ but we will only care about the integer $k(U)$ here). Dually, any quotient of $R_n$ is isomorphic to $I \oplus R_{k'(U)}$ for some preinjective $I$ and $0 \leq k'(U) \leq n$.

Cerulli Irelli and Esposito show that the loci
\[
X_d = \setst{U \in \Gr{R_n}{\beta}{}}{k(U), k'(U) \geq d}
\]
stratify $\Gr{R_n}{\beta}{}$ and that each stratum $X_d \setminus X_{d-1}$ is isomorphic to a classical Grassmannian variety (and thus has a cellular decomposition) \cite{Irelli:2010fk}.  We will show in this example how these strata can be constructed as rank loci.

The preprojective $P_n$ is the string module associated to the quiver mapping to $Q$:
\begin{equation}\label{eq:string}
\comment{
\begin{tikzpicture}
\draw (6,0) node {$2$}; 
\draw (5,1) node {$1$}; 
\draw (4,0) node {$2$}; 
\draw (3,1) node {$1$}; 
\draw (2,0) node {$2$}; 
\draw (1,1) node {$1$}; 
\draw (0,0) node {$2$}; 
\node (5) at (5,1) {$\bullet$};
\node (4) at (4,0) {};
\draw[->] (5) -- (4);
\end{tikzpicture}, }
2 \xleftarrow{b} 1 \xto{a} 2 \xleftarrow{b} 1 \xto{a} 2 \xleftarrow{b} 1 \xto{a} \cdots \xleftarrow{b} 1 \xto{a} 2
\end{equation}
(here, the labels of the vertices and arrows indicate where they map in $Q$).  More precisely, if we write write $p_n$ for this map of quivers when the string (\ref{eq:string}) has $n$ vertices marked 1,
and denote by $\id$ the representation of (\ref{eq:string}) with the vector space $K$ at every vertex and identity map over each arrow, then we get $P_n = p_{n*}(\id)$ using the pushforward definition from (\ref{eq:pf}).

By removing the first vertex marked 2 and adjacent arrow marked $b$, we get the string associated to the regular module $R_n$; denote the corresponding morphism from the string to $Q$ by $r_n$. This gives us two rank functions
\[
\rank_{p_n}, \rank_{r_n} \colon \repq \to \Z_{\geq 0} 
\]
on $Q$, from the definition (\ref{eq:pbrank}).
One can calculate the values of these rank functions on the representations $P_n$ and $R_n$ to be:
\begin{align}
\rank_{r_d}(R_n) = \maxint{n - d +1}  \qquad &\rank_{r_d}(P_n) = \maxint{n - d +1}  \\
\rank_{p_d}(R_n) = \maxint{n - d}  \qquad &\rank_{p_d}(P_n) = \maxint{n - d +1} .
\end{align}
From this we see that for $U \simeq P \oplus R_k$, we have
\[
\rank_{r_d} (U) - \rank_{p_d}(U) =
\begin{cases}
1 & k \geq d \\
0 & k <d \\
\end{cases},
\]
so that $k(U) \geq d$ if and only if $\rank_{r_d} (U) = \rank_{p_d}(U) +1$.  Since $Q^{\rm op} = Q$, it is easy to see by duality that $k'(U) \geq d$ if and only if $\rank_{r_d} (V/U) = \rank_{i_d}(V/U) + 1$, where $i_n$ is the morphism of quivers giving the preinjective $I_n$ as a string module.
So the stratification given by Cerulli Irelli and Esposito can be described by the rank loci
\begin{equation}
X_d = \setst{U \in \Gr{R_n}{\beta}{}}{\rank_{r_d} (U) = \rank_{p_d}(U) +1 = \rank_{i_d}(V/U) + 1} .
\end{equation}
%\todo{can we get the cells even as rank loci???}
\end{example}

\section{Future Directions}

Many of the natural questions suggested by the main result and examples fall under the general umbrella of: ``How do we choose rank functions and inequalities on them to get rank loci which are interesting in various situations?''  
In Example \ref{eg:strat}, more specifically we would like rank loci in a quiver Grassmannian which are better behaved or better understood than the original variety. In this example, we saw that a certain choice of rank data gave a stratification with strata isomorphic to known varieties (classical Grassmannians); more generally we might hope to construct rank loci which are at least fibered over some rank loci in a representation space of smaller dimension vector or for a smaller quiver.  One end goal would be computation of or positivity of Euler characteristics for quiver Grassmannians relevant to cluster algebras.

Example \ref{eg:typea} suggests a similar line of approach to the study of orbit closures in $\rep(\alpha)$.  Typically, there are  infinitely many orbits in a representation space if $Q$ is not of Dynkin or affine Dynkin type, and we currently have no clear picture of the orbits, much less how their closures relate.  Rank loci agglomerate many orbits by fixing discrete data; perhaps for wild-type quivers we can choose rank data to get loci whose degeneration order is more manageable.

Finally, it might be interesting to see if the singularities of certain rank loci (and their closures) are better behaved than orbit closures in wild representation type (see \cite{Zorbitclosure,Corbitclosures} for examples of bad singularities in infinite type, and the numerous papers by Zwara and Zwara-Bobi\'nski on singularities of orbit closures in more generality, e.g. \cite{MR1949357, MR1967381}).

\comment{
\section{Semi-continuous functions from rank functions}

\todo{probably just cut these last two non-sections}

\todo{When are rank functions sc? (only when actually rank) how to build sc functions from rank functions (can explicitly do hom on rooted trees, when others?)}

\section{Application?}

\todo{Think about example \cite[9.2]{Lusztig:1991yq}}

\todo{ the varieties $E_{=\beta}$ are conical in the vector space $\rep(\alpha)$ so we can projectivize them... used by caldero=Keller in ext groups}
}

\bibliographystyle{alpha}
\bibliography{ryanbiblio}

\def\cprime{$'$} \def\ocirc#1{\ifmmode\setbox0=\hbox{$#1$}\dimen0=\ht0
  \advance\dimen0 by1pt\rlap{\hbox to\wd0{\hss\raise\dimen0
  \hbox{\hskip.2em$\scriptscriptstyle\circ$}\hss}}#1\else {\accent"17 #1}\fi}
\begin{thebibliography}{DWZ09}

\bibitem[ASS06]{assemetal}
Ibrahim Assem, Daniel Simson, and Andrzej Skowro{\'n}ski.
\newblock {\em Elements of the representation theory of associative algebras.
  {V}ol. 1}, volume~65 of {\em London Mathematical Society Student Texts}.
\newblock Cambridge University Press, Cambridge, 2006.
\newblock Techniques of representation theory.

\bibitem[BZ02]{MR1967381}
Grzegorz Bobi{\'n}ski and Grzegorz Zwara.
\newblock Schubert varieties and representations of {D}ynkin quivers.
\newblock {\em Colloq. Math.}, 94(2):285--309, 2002.

\bibitem[Chi07]{Corbitclosures}
Calin Chindris.
\newblock On orbit closures for infinite type quivers.
\newblock \href{http://arXiv.org/abs/0709.3613}{\texttt{arxiv:0709.3613}},
  2007.

\bibitem[DWZ09]{Derksen:2009kl}
Harm Derksen, Jerzy Weyman, and Andrei Zelevinsky.
\newblock Quivers with potentials and their representations {II}:
  {A}pplications to cluster algebras.
\newblock \href{http://arXiv.org/abs/0904.0676}{\texttt{arxiv:0904.0676}},
  2009.

\bibitem[Ful93]{MR1234037}
William Fulton.
\newblock {\em Introduction to toric varieties}, volume 131 of {\em Annals of
  Mathematics Studies}.
\newblock Princeton University Press, Princeton, NJ, 1993.
\newblock The William H. Roever Lectures in Geometry.

\bibitem[Har77]{MR0463157}
Robin Hartshorne.
\newblock {\em Algebraic geometry}.
\newblock Springer-Verlag, New York, 1977.
\newblock Graduate Texts in Mathematics, No. 52.

\bibitem[IE10]{Irelli:2010fk}
Giovanni~Cerulli Irelli and Francesco Esposito.
\newblock Geometry of quiver grassmannians of kronecker type and canonical
  basis of cluster algebras.
\newblock \href{http://arxiv.org/abs/1003.3037}{\texttt{arXiv:1003.3037}},
  2010.

\bibitem[Kac83]{MR718127}
Victor~G. Kac.
\newblock Root systems, representations of quivers and invariant theory.
\newblock In {\em Invariant theory ({M}ontecatini, 1982)}, volume 996 of {\em
  Lecture Notes in Math.}, pages 74--108. Springer, Berlin, 1983.

\bibitem[Kel08]{Kclusterquivertriang}
Bernhard Keller.
\newblock Cluster algebras, quiver representations and triangulated categories.
\newblock \href{http://arXiv.org/abs/0807.1960}{\texttt{arxiv:0807.1960}},
  2008.

\bibitem[Kin08]{kinserrank}
Ryan Kinser.
\newblock The rank of a quiver representation.
\newblock {\em J. Algebra}, 320(6):2363--2387, 2008.

\bibitem[Kin09]{Kinser:2009zr}
Ryan Kinser.
\newblock {\em Rank functors and representation rings of quivers}.
\newblock PhD thesis, University of Michigan, 2009.
\newblock Available at author's webpage, currently:
  \url{http://www.math.uconn.edu/~kinser/thesis-single.pdf}.

\bibitem[Kin10]{kinserrootedtrees}
Ryan Kinser.
\newblock Rank functions on rooted tree quivers.
\newblock {\em Duke Math. J.}, 152(1):27--92, 2010.

\bibitem[KR86]{MR897322}
H.~Kraft and Ch. Riedtmann.
\newblock Geometry of representations of quivers.
\newblock In {\em Representations of algebras ({D}urham, 1985)}, volume 116 of
  {\em London Math. Soc. Lecture Note Ser.}, pages 109--145. Cambridge Univ.
  Press, Cambridge, 1986.

\bibitem[LeB88]{LeBruyn1986d}
Lieven LeBruyn.
\newblock Centers of generic division algebras and zeta-functions.
\newblock volume 1328 of {\em Lecture Notes in Mathematics}, 1988.
\newblock UIA report 86-24.

\bibitem[Nak96]{Nakajima:1996ys}
Hiraku Nakajima.
\newblock Varieties associated with quivers.
\newblock In {\em Representation theory of algebras and related topics
  ({M}exico {C}ity, 1994)}, volume~19 of {\em CMS Conf. Proc.}, pages 139--157.
  Amer. Math. Soc., Providence, RI, 1996.

\bibitem[Rei08]{Reineke:2008fk}
Markus Reineke.
\newblock Moduli of representations of quivers.
\newblock \href{http://arxiv.org/abs/0802.2147}{\texttt{arxiv:0802.2147}},
  2008.

\bibitem[Sch92]{schofieldgeneralreps}
Aidan Schofield.
\newblock General representations of quivers.
\newblock {\em Proc. London Math. Soc. (3)}, 65(1):46--64, 1992.

\bibitem[Zwa02]{MR1949357}
Grzegorz Zwara.
\newblock Unibranch orbit closures in module varieties.
\newblock {\em Ann. Sci. \'Ecole Norm. Sup. (4)}, 35(6):877--895, 2002.

\bibitem[Zwa03]{Zorbitclosure}
Grzegorz Zwara.
\newblock An orbit closure for a representation of the {K}ronecker quiver with
  bad singularities.
\newblock {\em Colloq. Math.}, 97(1):81--86, 2003.

\end{thebibliography}

\end{document}